\documentclass[journal]{new-aiaa}
\usepackage[utf8]{inputenc}

\usepackage{graphicx}
\usepackage{amsmath}
\usepackage[version=4]{mhchem}
\usepackage{siunitx}
\usepackage{longtable,tabularx}
\usepackage{hyperref}

\title{On the Efficiency of the Iterative Techniques for Solving Incompressible Navier-Stokes Equations}

\author{Mohamed Mohsen Ahmed\footnote{Corresponding author email: mmahmed@umd.edu}} 
\affil{Department of Mechanical Engineering, University of Maryland, College Park, Maryland, USA}
\begin{document}

\maketitle

\begin{abstract}
It is well known that the choice of the iterative method is crucial in determining the speed of the converged solution. This article presents a detailed comparison between several iterative techniques for solving incmopressible Navier-Stokes equations. The numerical approaches implemented in the solver include Jacobi, Gauss-Siedel, Successive Over Relaxation, Alternating Direct Implicit and Multigrid methods. The results reveal that multigrid method is the most powerful iterative method among all other methods investigated in terms of the computational time and the number of iterations. 
\end{abstract}


\section{Introduction}
Navier-Stokes equations are non-linear second order coupled partial differential equations that an exact solution can only be obtained if certain approximations are performed. Therefore, numerical solutions of Navier-Stokes equations are the most practical method for engineering application involving complex flows. However, complex flows are associated with different length and time scales such that a large number of iterations using a sufficiently large number of grid points or mesh cells should be carried out in order to converge to a solution that captures the complex flow physics. Although parallel computing has been widely used to deal with large number of mesh cells by distributing the mesh on several processors using different decomposition methods \cite{shang2014impact}, the choice of the iterative method that shall be used to solve the governing equations remains crucial in determining the speed of obtain a converged solution.

This paper presents a detailed comparison of the speed of the most common iterative techniques that are employed in solving a system of equations. The iterative techniques employed in this study include Jacobi, Gauss-Siedel (GS), Successive Over Relaxation (SOR), Successive Over Relaxation by Line (SLOR), Alternative Direct Implicit (ADI) and multigrid. The numerical example employed in this study is the solution of the incompressible Navier-Stokes equations using projection method \cite{brown2001accurate}. 

The governing equations in the non-dimensionalized primitive variable form are

\begin{enumerate}
\item continuity equation
\begin{equation}
\frac{\partial u}{\partial x}+\frac{\partial v}{\partial y}=0
\end{equation}

\item Momentum equation x-direction
\begin{equation}
\frac{\partial u}{\partial t}+ \frac{\partial u^2}{\partial x} + \frac{\partial uv}{\partial y} + \frac{\partial p}{\partial x}= \frac{1}{Re} (\frac{\partial^2 u} {\partial x^2} + \frac{\partial^2 u} {\partial y2})
\end{equation}

\item Momentum equation y-direction
\begin{equation}
\frac{\partial v}{\partial t}+ \frac{\partial uv}{\partial x} + \frac{\partial v^2} {\partial y} + \frac{\partial p}{\partial y}= \frac{1}{Re} (\frac{\partial^2 v} {\partial x^2} + \frac{\partial^2 v} {\partial y2})
\end{equation}

\end{enumerate}

\section{Discretization of the Governing Equations}
Finite difference method is employed to discretize  the governing equations on a staggered grid. The use of staggered grid allows the coupling between $u,v$ and $p$ at adjacent point, thus, preventing oscillations and checker-board effects. In the description of staggered grid, the pressure is defined at the cell center point $i,j$. However, the velocity components are defined as $u_f, u_b , v_f$ and $v_b$ which correspond to forward and backward edges of the cell corresponding to point $i,j$. The velocity at the cell center point $i,j$ is calculated from the average of the forward and backward velocities such that
\begin{align}
u(i,j)=\frac{u_f(i,j)+u_b(i,j)}{\Delta x}\\
v(i,j)=\frac{v_f(i,j)+v_b(i,j)}{\Delta y}
\end{align}

The continuity equation can be discretized using central differencing about point $i,j$ as
\begin{equation}
\frac{u_f(i,j)-u_b(i,j)}{\Delta x} + \frac{v_f(i,j)-v_b(i,j)}{\Delta y}=0
\end{equation}

Explicit discretization of the x-momentum equation about point $i,j$ using central differences yields
\begin{align}
u^{k+1}_f(i,j)=F_f^{k} (i,j)-\frac{\Delta t}{\Delta x}  \big[p(i+1,j)-p(i,j)\big]^{k+1}
\end{align}
where
\begin{align}
F_f^k(i,j) =u_f^k(i,j) &+ \Delta t \Big[\frac{u_f(i+1,j)-2u_f(i,j)+u_b(i,j)}{Re \Delta x^2} \nonumber\\
&+ \frac{u_f(i,j-1)-2u_f(i,j)+u_f(i,j+1)}{Re \Delta y^2} \\
&-\frac{u^2(i+1,j)-u^2(i,j)}{\Delta x} - \frac{[uv]_{ff}(i,j)-[uv]_{fb}(i,j)}{\Delta y} \Big]^k \nonumber
\end{align}

Similarly, the y-momentum equation is discretized as
\begin{align}
v^{k+1}_f(i,j)=G_f^{k} (i,j)-\frac{\Delta t}{\Delta x}  \big[p(i,j+1)-p(i,j)\big]^{k+1}
\end{align}
where
\begin{align}
G_f^k(i,j) =v_f^k(i,j) &+ \Delta t \Big[\frac{v_f(i+1,j)-2v_f(i,j)+v_f(i-1,j)}{Re \Delta x^2} \nonumber \\
&+ \frac{u_f(i,j+1)-2v_f(i,j)+v_b(i,j)}{Re \Delta y^2} \\
&-\frac{v^2(i,j+1)-v^2(i,j)}{\Delta y} - \frac{[uv]_{ff}(i,j)-[uv]_{bf}(i,j)}{\Delta x} \Big]^k \nonumber
\end{align}
The subscripts $ff,fb$ and $bf$ correspond to the corner points of the cell whose center is point $i,j$. The corresponding terms are
\begin{align}
[uv]_{ff}(i,j)&=\big(\frac{u_f(i,j)-u_f(i,j+1)}{2} \big) \big(  \frac{v_f(i,j)+v_f(i+1,j)}{2} \big) \nonumber \\
[uv]_{fb}(i,j)&=\big(\frac{u_f(i,j)-u_f(i,j-1)}{2} \big) \big(  \frac{v_f(i,j-1)+v_f(i+1,j-1)}{2} \big) \\
[uv]_{bf}(i,j)&=\big(\frac{u_f(i-1,j)-u_f(i-1,j+1)}{2} \big) \big(  \frac{v_f(i,j)+v_f(i-1,j)}{2} \big) \nonumber
\end{align}

The pressure Poisson equation is finally obtained from the discretized continuity and momentum equations such that
\begin{equation}
\Big[ \frac{p(i-1,j)-2p(i,j) +p(i+1,j) }{\Delta x^2} +  \frac{p(i,j-1)-2p(i,j) +p(i,j+1) }{\Delta y^2}     \Big]^{k+1} = RHS^k(i,j)
\end{equation}
where
\begin{equation}
RHS^k(i,j)=\frac{1}{\Delta t} \Big[ \frac{F_f(i,j)-F_b(i,j)}{\Delta x} + \frac{G_f(i,j)-G_b(i,j)}{\Delta y} \Big]^k
\end{equation}
Since projection method provides an explicit expression of the momentum equation, limitations are introduced on the maximum time step for a stable solution. Peyret and Taylor \cite{peyret} provided a correlation between the time step, Reynolds number and mesh size. If $\Delta x=\Delta y$ the restriction is
\begin{align}
\frac{\Delta t}{Re \ \Delta x^2} \le 0.25
\end{align}

\begin{figure}[h]
\centering
\includegraphics[width=0.5\textwidth]{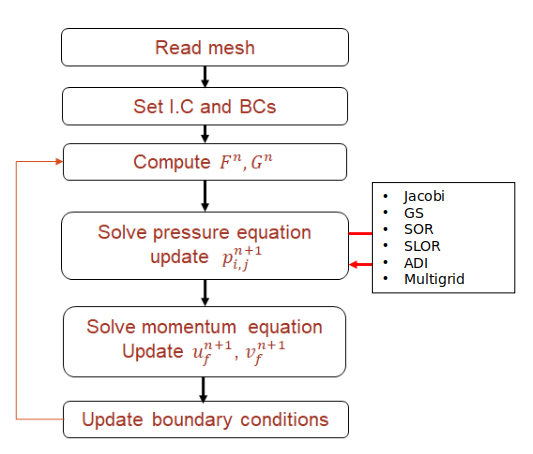} 
\caption{Solver algorithm} \label{fig3}
\end{figure}

Now we have all equations in a discretized form, the next step is to implement these equations in a FORTRAN code. The overall FORTRAN solver algorithm is explained in the flowchart of figure \ref{fig3}. 

\section{Results and Discussion}
\begin{figure}[h!]
\centering
\begin{minipage}{.65\textwidth}
  \centering
  \includegraphics[width=1.0\linewidth]{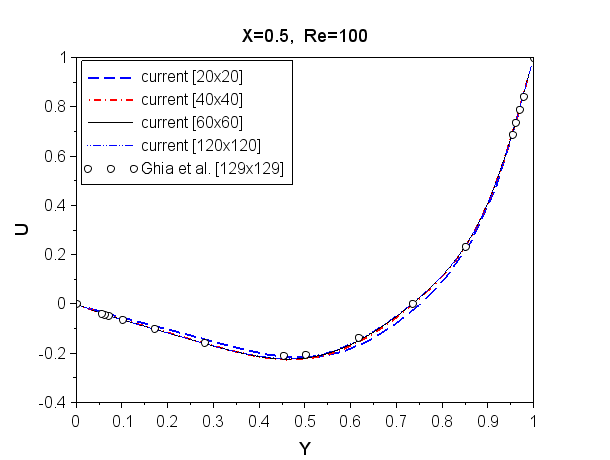}
  \captionof{figure}{distribution of $u$ at $x=0.5$ for $Re=100$}
  \label{fig7}
\end{minipage}%
\begin{minipage}{.35\textwidth}
  \centering
  \includegraphics[width=0.9\linewidth]{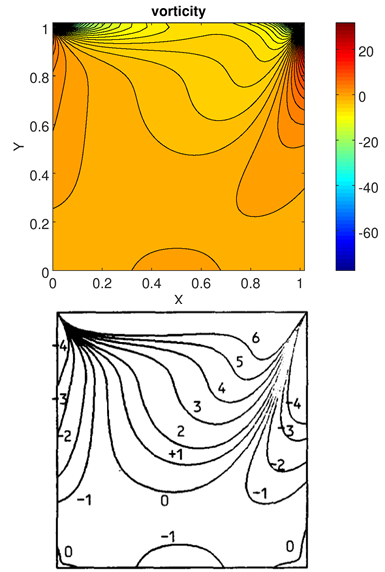}
  \captionof{figure}{Re 100 Vorticity contours, top: current solution, bottom Ghia et al. \cite{ghia}}
  \label{fig8}
\end{minipage}
\end{figure}

Two test cases are considered. First, a square cavity of unity length with a moving top surface is considered for validation and verification purposes. The top wall is moving with speed $v_w=1$ and all other walls are not moving. Different grid sizes are considered for each simulated case. Ghost cells are introduced on the staggered grid to specify the constraints of the boundary points. No slip and zero pressure gradient are specified on all walls. A study on the mesh size effect on the solution is performed first at a Reynolds number of 100. A qualitative comparison is carried out between the current solution and two different benchmark solutions by Ghia et al. \cite{ghia}. Figure \ref{fig7} shows the distribution of $u$ velocity along y-direction at $x=0.5$ for different mesh levels at $Re=100$ when the solution reached steady state (time=17.24 s). It can be inferred that a grid independent solution can be obtained at mesh level 60x60 or above. 

\begin{figure}[h]
\centering
  \centering
  \includegraphics[width=.4\linewidth]{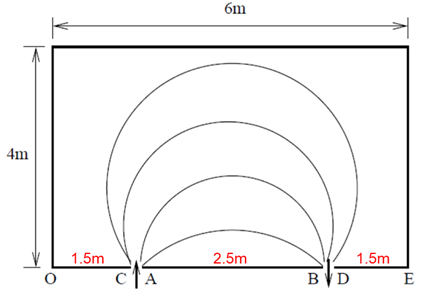}
  \captionof{figure}{Geometry description}
  \label{case2}
\end{figure}

The second test case involves flow in the chamber described in figure ~\ref{case2}. The inlet and outlet small gaps are $CA=BD=0.25m$. A uniform grid is constructed with constant spacing $\Delta x = \Delta y=0.25 m$. The computational grid and boundary condition is provided in figure ~\ref{fig3}. The total number of grid pints in the domain is 425 points. The simulations are performed on a single processor and the CPU time is recorded for all cases.

\begin{figure}[h]
\centering
\begin{minipage}{.5\textwidth}
  \centering
  \includegraphics[width=.9\linewidth]{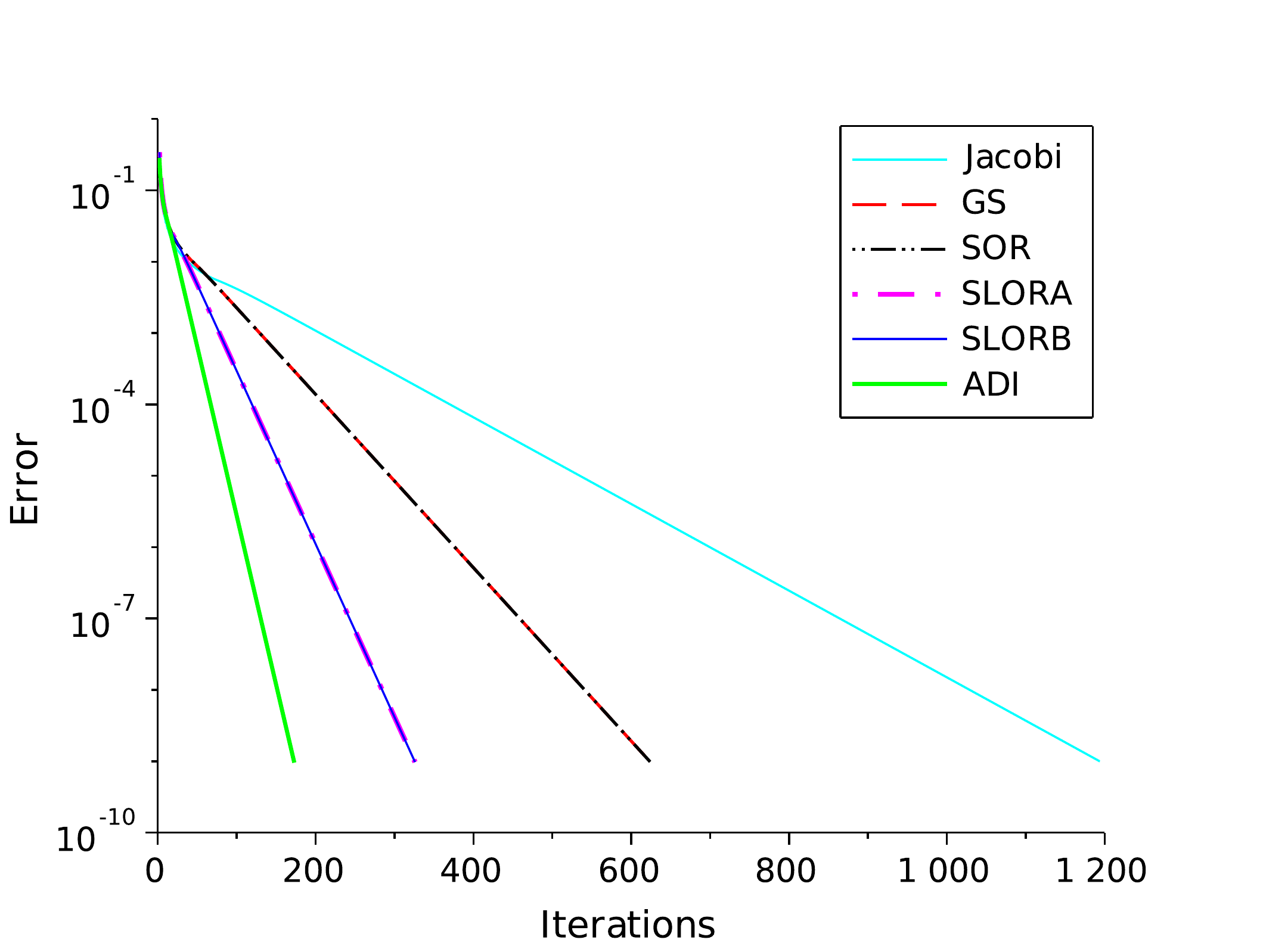}
  \captionof{figure}{Error vs. iterations (no over relaxation)}
  \label{fig11}
\end{minipage}%
\begin{minipage}{.5\textwidth}
  \centering
  \includegraphics[width=1.0\linewidth]{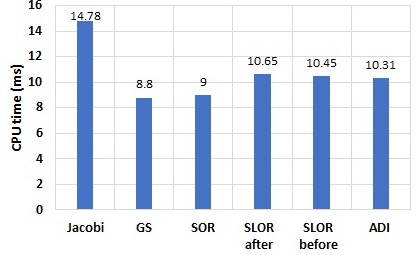}
  \captionof{figure}{Computational time (no over relaxation)}
  \label{fig12}
\end{minipage}
\end{figure}

All numerical methods are first compared without over relaxation in order to give us confidence on the validity of our solver to predict well known trends. Figure ~\ref{fig11} shows the number of iterations and the error of the different iterative methods employed in the solver. The slowest method is Jacobi which requires 1194 iterations to reach convergence. GS and SOR with no relaxation are identical and the number of iterations is 624 iterations, which is approximately half the number of iterations of Jacobi method. SLORB and SLORA, which refer to applying over relaxation before or after each GS line iteration cycle, are much faster than GS and SOR method with only 326 iterations. Although ADI method is the fastest among these methods with only 326 iterations, it should be taken into consideration that two calculations are conducted inside a single cycle of the ADI method. In figure ~\ref{fig12}, the CPU time for each method is reported. It can be inferred that Jacobi is the slowest method, GS is the fastest method and SLOR and ADI methods are approximately equal in terms of computational time due to the reasons discussed previously. 

\begin{figure}[h]
\centering
\begin{minipage}{.5\textwidth}
  \centering
  \includegraphics[width=1\linewidth]{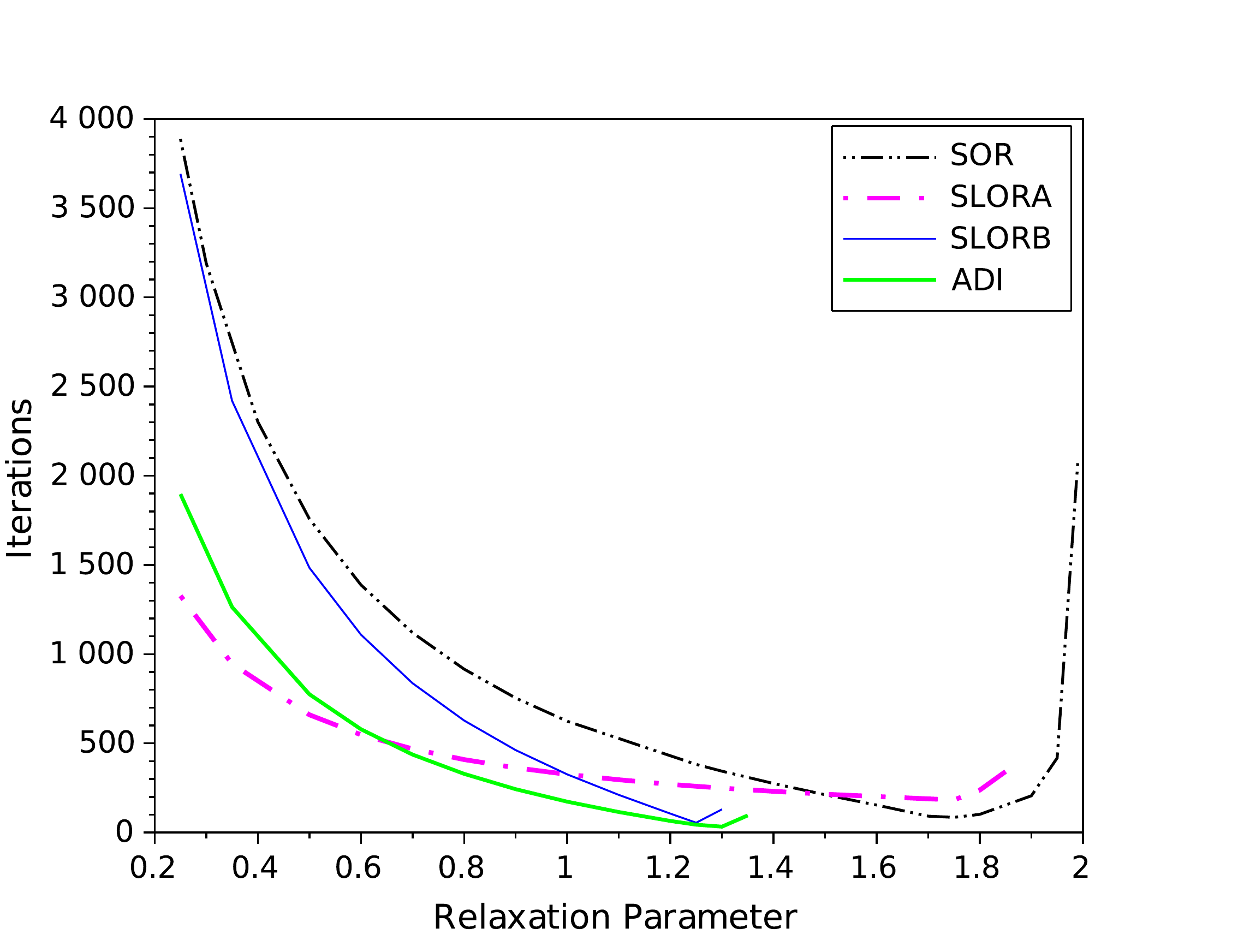}
  \captionof{figure}{iterations vs. relaxation parameter $w$}
  \label{fig13}
\end{minipage}%
\begin{minipage}{.5\textwidth}
  \centering
  \includegraphics[width=1\linewidth]{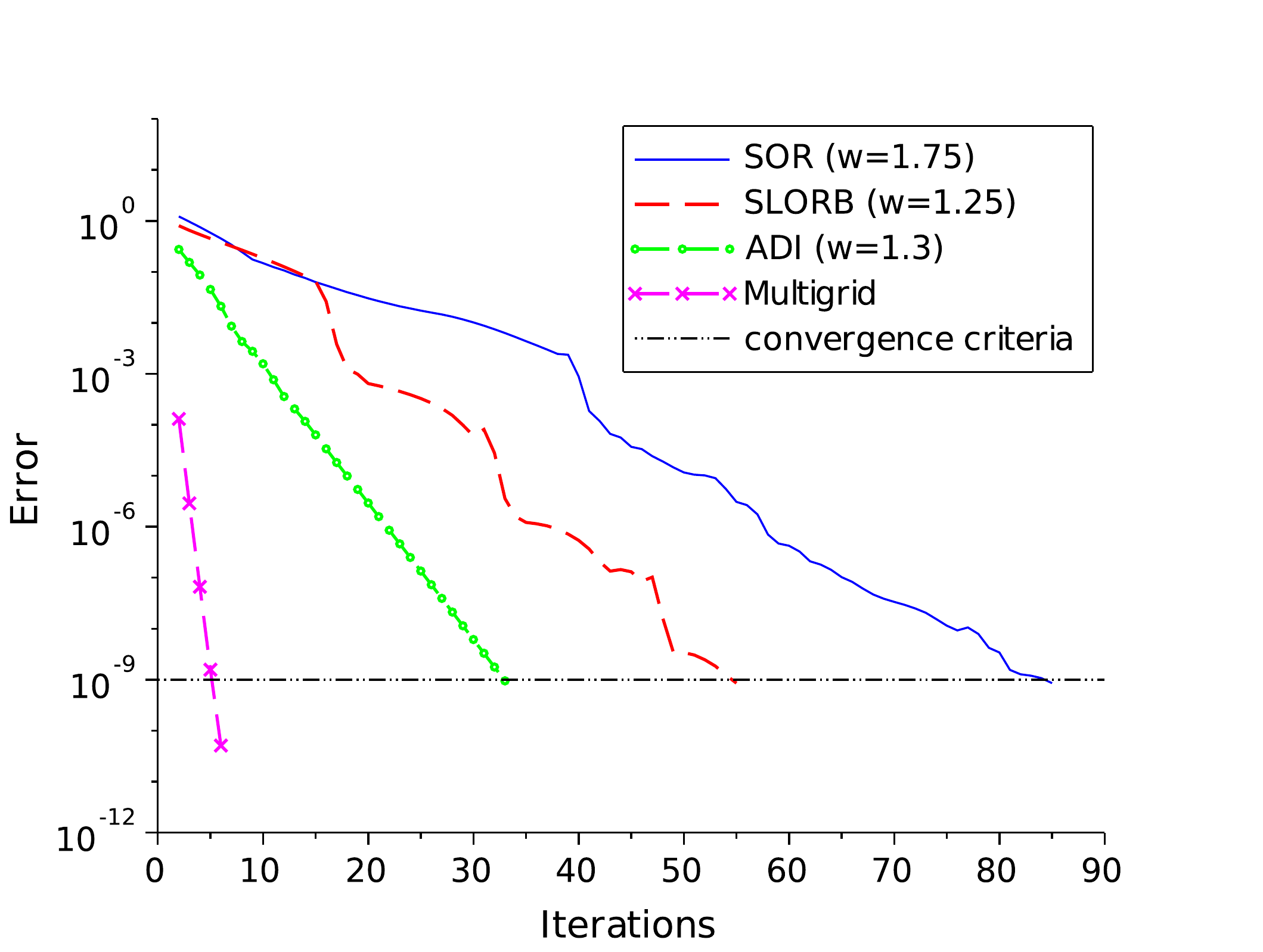}
  \captionof{figure}{Optimum relaxed methods}
  \label{fig14}
\end{minipage}
\end{figure}

Now, we are in turn to find the value of the best over relaxation parameter that would speed up the calculations of the SOR, SLOR and ADI methods. Figure ~\ref{fig13} show the effect of the relaxation parameter on the number of iterations. It can be seen that there exist an optimum value of $1<w<2$ and this value is different for different iteration methods. The optimum values of $w$ are 1.75, 1.75, 1.25 and 1.3 for SOR, SLORA, SLORB and ADI respectively. Figure ~\ref{fig14} shows the error vs. number of iterations required at the best relaxation parameter $w$ for each iterative method. \\

\begin{figure}[h]
\centering
\includegraphics[width=.8\textwidth]{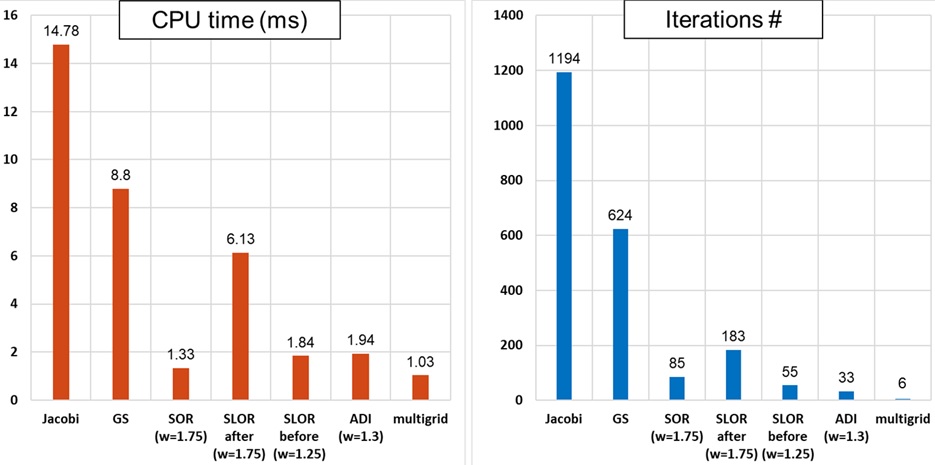} 
\caption{Summary of Iteration number and CPU time} \label{fig15}
\end{figure}

\begin{figure}[h]
\centering
\includegraphics[width=.5\textwidth]{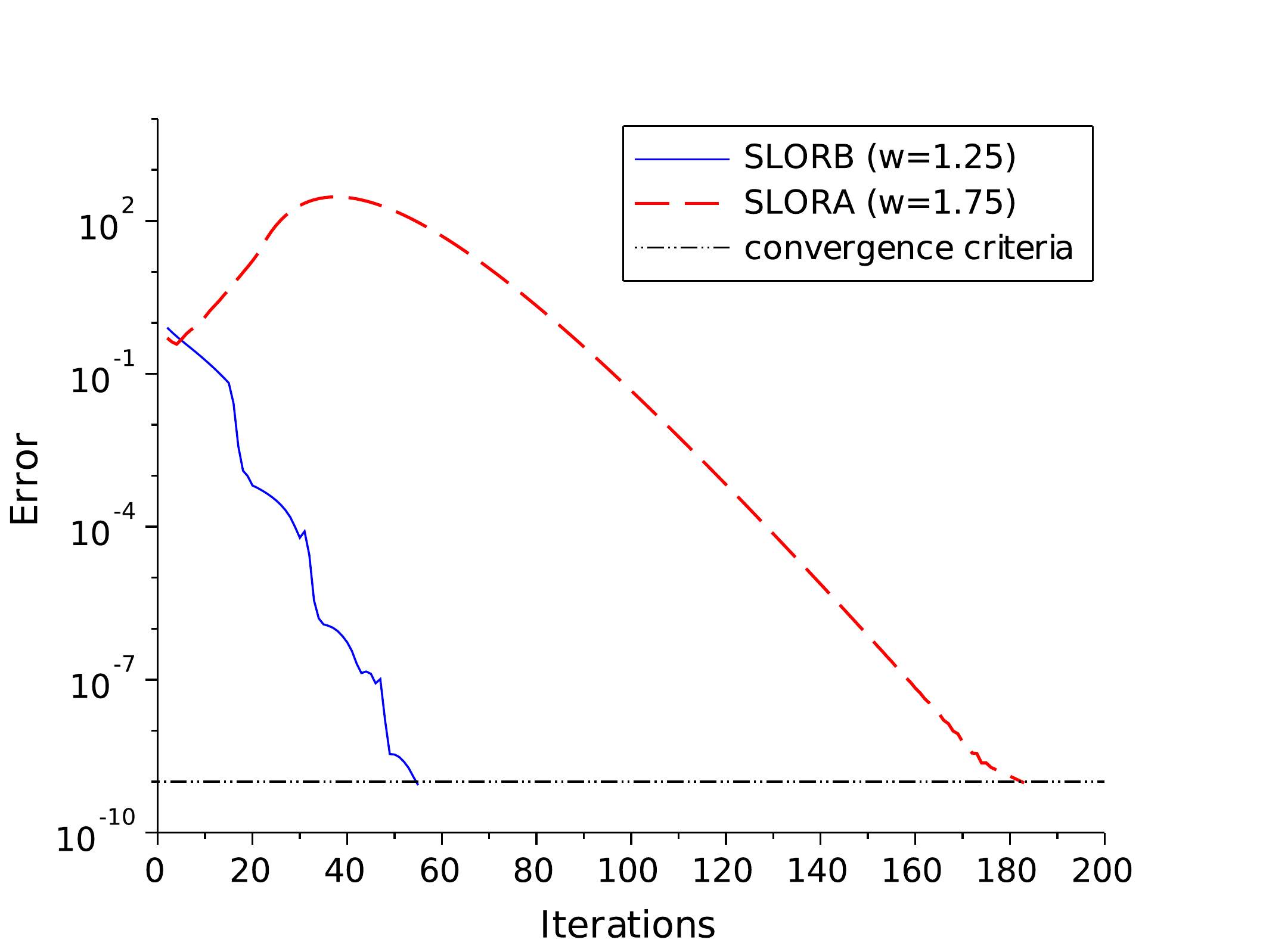} 
\caption{SLORB and SLORA with optimal relaxation parameter} \label{fig16}
\end{figure}

The number of iterations and CPU times for all methods with optimal relaxation are summarized in figure ~\ref{fig15}. As expected, ADI requires less number of iterations that SLORB (33 vs. 55). However, the total CPU time is almost identical in both cases. Another important observation is that SLORA is much slower than SLORB and SOR with optimal over relaxation parameter. Figure ~\ref{fig16} compares SLORB and SLORA at their optimal over relaxation parameters. It is inferred that although $w$ is higher in SLORA, more number of iterations is required since the solution overshoots in the first 50 iterations. The only explanation of this observation is low stability of Thomas' algorithm \cite{thomas1949elliptic} while solving the tri-diagonal matrix when the matrix is not diagonally dominant.

It is also concluded from figure ~\ref{fig15} that multigrid method is the most efficient method among all other iterative methods. The solution converged after only six V-cycles and the total CPU time is around 1.03 milliseconds which is fastest relative to all other methods. 

\section{Conclusion}
Despite the rapid growth in parallel computing, the properties of the iterative methods are the key parameter for determining the speed of obtaining a converged solution on a single processor. Six different implicit and explicit iterative methods were compared while solving a 2-dimensional incompressible Navier-Stokes eqations. The results reveal multigrid method is more efficient that other iterative methods in terms of number of iterations and computational time. 

\section{Funding}
This study was conducted without a funding support.

\section{Conflict of interest}
The authors declare that there is no conflict of interest.

\newpage
\section{Appendix: Incompressible N-S Equations Fortran Solver}
\begin{verbatim}
 program Incomp_NS

implicit none
!!!!!!-------variables declaration--------!!!!!
   integer                            :: m,n,i,j, iteration , counter ,anim_freq
   real(8)                            :: Lx,Ly,dx,dy,beta, v_w , Ren, cycles ,relaxation
   real(4)                            :: time, t_final, dt , CFL
   real(8),allocatable,dimension(:,:) :: P, u_f, u_b, v_f , v_b , u , v , R , psi , vor
   real(8),allocatable,dimension(:,:) :: F_f, F_b , G_f, G_b
   real(8),allocatable,dimension(:,:) :: uv_ff , uv_fb , uv_bf , uv_bb
   real(8),allocatable,dimension(:,:) :: P_old,P_new, difference
   real(8)                            :: max_error , u_center
   character(len=50)                  :: filename
   real(8)                            :: mm,nn
!!!!!-------------Inputs-------------!!!!!
cycles=100000
anim_freq=500
Lx=1.0d0
Ly=1.0d0
v_w=1.0d0
m=240
n=240
dx=Lx/(m-2)
dy=Ly/(n-2)
beta=dx/dy
time=0.0d0
Ren=10000.0d0
dt=0.0025*Ren*dx**2
t_final=cycles*dt
CFL=dt/dx
relaxation=1.3
mm=m
nn=n

!!!!!-------Creating variables--------!!!!!
allocate(P(m,n) , u(m,n) , v(m,n) , R(m,n) , psi(m,n) , vor(m,n))
allocate(u_f(m,n) , u_b(m,n) , v_f(m,n) , v_b(m,n))
allocate(uv_ff(m,n) , uv_fb(m,n) , uv_bf(m,n) , uv_bb(m,n))
allocate(F_f(m,n), F_b(m,n), G_f(m,n), G_b(m,n))
allocate(P_new(m,n),P_old(m,n),difference(m,n))
P=0.0d0
u=0.0d0
v=0.0d0
R=0.0d0
psi=0.0d0
vor=0.0d0
uv_ff=0.0d0
uv_fb=0.0d0
uv_bf=0.0d0
uv_bb=0.0d0
u_f=0.0d0
u_b=0.0d0
v_f=0.0d0
v_b=0.0d0
F_f=0.0d0
F_b=0.0d0
G_f=0.0d0
G_b=0.0d0
P_old=0.0d0
P_new=0.0d0
!!!!!-------Update Boundary conditions--------!!!!!
    do i=1,m
       !bottom
       v_f(i,1)=0.0d0
       u_f(i,1)= -u_f(i,2)
       P(i,1)=P(i,2)
       !top
       v_f(i,n-1)=0.0d0
       u_f(i,n)=2*v_w-u_f(i,n-1)
       P(i,n)=P(i,n-1)
    end do
    do j=1,n
        !left
        u_f(1,j)=0.0d0
        v_f(1,j)= -v_f(2,j)
        P(1,j)=P(2,j)

        !right
        u_f(m-1,j)=0.0d0
        v_f(m,j)= -v_f(m-1,j)
        P(m,j)=P(m-1,j)
    end do
    do i=2,m
        do j=1,n
            u_b(i,j)=u_f(i-1,j)
        end do
    end do
    do i=1,m
        do j=2,n
            v_b(i,j)=v_f(i,j-1)
        end do
    end do
    do i=1,m
        do j=1,n
            u(i,j)=0.5*(u_f(i,j)+u_b(i,j))
            v(i,j)=0.5*(v_f(i,j)+v_b(i,j))
        end do
    end do

!!!!!----------------------------Main Solver------------------------------!!!!!
counter=0
      write (filename,'(A,I6.6,A)') "p", counter,".csv"
      print*,filename

      open(unit=100,file=filename)
      do i=2,m-1
         write(100,*) (P(i,j),',', j = 2,n-2), P(i,n-1)
      enddo
      close(unit=100)

      write (filename,'(A,I6.6,A)') "u", counter,".csv"
      print*,filename
      open(unit=200,file=filename)
      do i=2,m-1
         write(200,*) (u(i,j),',', j = 2,n-2), u(i,n-1)
      enddo
      close(unit=200)

      write (filename,'(A,I6.6,A)') "v", counter,".csv"
      print*,filename
      open(unit=300,file=filename)
      do i=2,m-1
         write(300,*) (v(i,j),',', j = 2,n-2), v(i,n-1)
      enddo
      close(unit=300)

open(unit=10,file="Time_U.csv")

do counter=1,cycles

    time=time+dt
    u_center=u(nint(mm/2),nint(nn/4))
    print*, '*********************************************'
    print*,'cycle number =    ', counter
    print*, 'dt= ',dt
    print*, 'dx= ',dx
    print*,'current time =    ', time
    print*,'Courant number =    ', CFL
    print*,'monitoring velocity=  ', u_center

    do i=2,m-1
        do j=2,n-1
            uv_ff(i,j)=0.5*(u_f(i,j)+u_f(i,j+1)) * 0.5*(v_f(i,j)+v_f(i+1,j))
            uv_fb(i,j)=0.5*(u_f(i,j)+u_f(i,j-1)) * 0.5*(v_b(i,j)+v_b(i+1,j))
            uv_bf(i,j)=0.5*(u_b(i,j)+u_b(i,j+1)) * 0.5*(v_f(i,j)+v_f(i-1,j))
            uv_bb(i,j)=0.5*(u_b(i,j-1)+u_b(i,j)) * 0.5*(v_b(i,j)+v_b(i-1,j))
        end do
    end do

! Intermediate step (F,G)
    do i=2,m-1
        do j=2,n-1
            F_f(i,j)=u_f(i,j) +(dt/(Ren*dx**2))*(u_f(i+1,j)-2*u_f(i,j)+u_b(i,j)) &
                     +(dt/(Ren*dy**2))*(u_f(i,j-1)-2*u_f(i,j)+u_f(i,j+1)) &
                     -(dt/dx)*(u(i+1,j)**2 - u(i,j)**2) &
                     -(dt/dy)*(uv_ff(i,j)-uv_fb(i,j))

            G_f(i,j)=v_f(i,j) +(dt/(Ren*dy**2))*(v_f(i+1,j)-2*v_f(i,j)+v_f(i-1,j)) &
                     +(dt/(Ren*dy**2))*(v_f(i,j+1)-2*v_f(i,j)+v_b(i,j)) &
                     -(dt/dx)*(uv_ff(i,j) - uv_bf(i,j)) &
                     -(dt/dy)*(v(i,j+1)**2 - v(i,j)**2)

            F_b(i,j)=u_b(i,j) +(dt/(Ren*dx**2))*(u_f(i,j)-2*u_b(i,j)+u_b(i-1,j)) &
                     +(dt/(Ren*dy**2))*(u_b(i,j-1)-2*u_b(i,j)+u_b(i,j+1)) &
                     -(dt/dx)*(u(i,j)**2 - u(i-1,j)**2) &
                     -(dt/dy)*(uv_bf(i,j)-uv_bb(i,j))

            G_b(i,j)=v_b(i,j) +(dt/(Ren*dy**2))*(v_b(i+1,j)-2*v_b(i,j)+v_b(i-1,j)) &
                     +(dt/(Ren*dy**2))*(v_f(i,j)-2*v_b(i,j)+v_b(i,j-1)) &
                     -(dt/dx)*(uv_fb(i,j) - uv_bb(i,j)) &
                     -(dt/dy)*(v(i,j)**2 - v(i,j-1)**2)

        end do
    end do


! Pressure equation
    do i=1,m
        do j=1,n
            R(i,j)=( (F_f(i,j)-F_b(i,j))/dx + (G_f(i,j)-G_b(i,j))/dy ) / dt
        end do
    end do

    !call an iterative solver
    call ADI(P,R,m,n,beta,dx,relaxation)

! Momentum equation
    do i=2,m-1
        do j=2,n-1
            u_f(i,j)=F_f(i,j) - dt/dx*(P(i+1,j) - P(i,j))
            v_f(i,j)=G_f(i,j) - dt/dy*(P(i,j+1) - P(i,j))
        end do
    end do


!!!!!-------Update Boundary conditions--------!!!!!
    do i=1,m
       !bottom
       v_f(i,1)=0.0d0
       u_f(i,1)= -u_f(i,2)
       P(i,1)=P(i,2)

       !top
       v_f(i,n-1)=0.0d0
       u_f(i,n)=2*v_w-u_f(i,n-1)
       P(i,n)=P(i,n-1)
    end do

    do j=1,n
        !left
        u_f(1,j)=0.0d0
        v_f(1,j)= -v_f(2,j)
        P(1,j)=P(2,j)

        !right
        u_f(m-1,j)=0.0d0
        v_f(m,j)= -v_f(m-1,j)
        P(m,j)=P(m-1,j)
    end do

    do i=2,m
        do j=1,n
            u_b(i,j)=u_f(i-1,j)
        end do
    end do

    do i=1,m
        do j=2,n
            v_b(i,j)=v_f(i,j-1)
        end do
    end do

    do i=1,m
        do j=1,n
            u(i,j)=0.5*(u_f(i,j)+u_b(i,j))
            v(i,j)=0.5*(v_f(i,j)+v_b(i,j))
        end do
    end do


!vorticity calculation

    do i=1,m-1
        do j=1,n-1
            psi(i,j)=psi(i-1,j)-dx*(v(i,j))
            vor(i,j)=(v(i+1,j)-v(i,j))/(dx) - (u(i,j+1)-u(i,j))/(dy)
        end do
    end do

!Animation output
   if(mod(counter,anim_freq)==0) then
      write (filename,'(A,I6.6,A)') "p", counter,".csv"
      print*,filename

      open(unit=100,file=filename)
      do i=2,m-1
         write(100,*) (P(i,j),',', j = 2,n-2), P(i,n-1)
      enddo
      close(unit=100)

      write (filename,'(A,I6.6,A)') "u", counter,".csv"
      print*,filename
      open(unit=200,file=filename)
      do i=2,m-1
         write(200,*) (u(i,j),',', j = 2,n-2), u(i,n-1)
      enddo
      close(unit=200)

      write (filename,'(A,I6.6,A)') "v", counter,".csv"
      print*,filename
      open(unit=300,file=filename)
      do i=2,m-1
         write(300,*) (v(i,j),',', j = 2,n-2), v(i,n-1)
      enddo
      close(unit=300)

   endif


write(10,*) time,',', u_center

end do

close(unit=10)

!!!!!-------Output results--------!!!!
open(unit=1, file="p.csv")
do i=2,m-1
write(1,*) (P(i,j),',', j = 2,n-2), P(i,n-1)
enddo
close(unit=1)

open(unit=2, file="u.csv")
do i=2,m-1
write(2,*) (u(i,j),',', j = 2,n-2), u(i,n-1)
enddo
close(unit=2)

open(unit=3, file="v.csv")
do i=2,m-1
write(3,*) (v(i,j),',', j = 2,n-2), v(i,n-1)
enddo
close(unit=3)

open(unit=4, file="stream.csv")
do i=2,m-1
write(4,*) (psi(i,j),',', j = 2,n-2), psi(i,n-1)
enddo
close(unit=4)

open(unit=5, file="vorticity.csv")
do i=2,m-1
write(5,*) (vor(i,j),',', j = 2,n-2), vor(i,n-1)
enddo
close(unit=5)

end program Incomp_NS
   
\end{verbatim}
\bibliography{sample}

\end{document}